\documentclass[11pt]{amsart}

\theoremstyle{plain}
\newtheorem{theorem}                {Theorem}      [section]
\newtheorem*{theorem*}                {Theorem \ref{thm:appl}}

\newtheorem{lemma}        [theorem]  {Lemma}

\theoremstyle{definition}

\newtheorem{remark}       [theorem]  {Remark}
\newtheorem{definition}   [theorem]  {Definition}

\DeclareMathOperator{\trace}{trace}

 \DeclareMathOperator{\id}{I}

\numberwithin{equation}{section}

\begin{document}

\title[Surfaces with parallel mean curvature]
{A note on surfaces with parallel mean curvature\vskip 1cm Une note sur des surfaces de courbure moyenne parall\`{e}le }

\author{Dorel~Fetcu}
\author{Harold~Rosenberg}

\address{Department of Mathematics\\
"Gh. Asachi" Technical University of Iasi\\
Bd. Carol I no. 11 \\
700506 Iasi, Romania} \email{dfetcu@math.tuiasi.ro}

\curraddr{IMPA\\ Estrada Dona Castorina\\ 110, 22460-320 Rio de
Janeiro, Brasil} \email{dorel@impa.br}

\address{IMPA\\ Estrada Dona Castorina\\ 110, 22460-320 Rio de
Janeiro, Brasil} \email{rosen@impa.ro}

\thanks{The first author was supported by a Post-Doctoral Fellowship "P\'os-Doutorado S\^enior
(PDS)" offered by FAPERJ, Brazil.}

\begin{abstract}
We use a Simons type equation in order to characterize complete non-minimal pmc surfaces with non-negative Gaussian curvature.

\vskip 0.5\baselineskip

Dans cette note, on \'{e}tudie des immersions isom\'etriques de surfaces compl\`{e}tes $\Sigma^2$ dans $M^n(c)\times\mathbb{R}$, ou $M^n(c)$ est une vari\' et\' e compl\`{e}te simplement connexe de courbure sectionnelle constante $c$. On classifie ces immersions, lorsque leur vecteur courbure moyenne est parall\`{e}le dans le fibr\'e normal et leur courbure intrins\`{e}que est positive ou nulle. L'outil principal est une diff\'erentielle quadratique holomorphe dont la partie sans trace satisfait l'\' equation de Codazzi.
\end{abstract}



\maketitle

\section{The main result}

Let $M^n(c)$ be a simply-connected $n$-dimensional manifold, with constant sectional curvature $c$, consider the product manifold $\bar M=M^n(c)\times\mathbb{R}$, and let $\Sigma^2$ be an immersed surface in $\bar M$.

\begin{definition} The surface $\Sigma^2$ is called a \textit{pmc surface} if its mean curvature vector $H$ is
parallel in the normal bundle. More precisely,
$\nabla^{\perp}H=0$, where the normal connection $\nabla^{\perp}$ is
defined by the Weingarten equation
$$
\bar\nabla_XV=-A_VX+\nabla^{\perp}_XV,
$$
for any vector field $X$ tangent to $\Sigma^2$ and any vector
field $V$ normal to the surface. Here $\bar\nabla$ is the
Levi-Civita connection on $\bar M$ and $A$ is the shape operator.
\end{definition}

When the dimension of $\bar M$ is equal to $3$, an immersed pmc surface
in $\bar M$ is a surface with constant mean curvature (a \textit{cmc surface}).
U. Abresch and H. Rosenberg introduced in \cite{AR,AR2} a holomorphic differential on such surfaces and then completely classified those cmc surfaces on which it vanishes.
In order to extend their results to the case of ambient spaces $\bar M=M^n(c)\times\mathbb{R}$,
with $n\geq 2$, H. Alencar, M. do Carmo and R. Tribuzy, defined in \cite{AdCT1, AdCT}
a real quadratic form $Q$ on pmc surfaces by
\begin{equation}\label{eq:Q}
Q(X,Y)=2\langle A_HX,Y\rangle-c\langle X,\xi\rangle\langle Y,\xi\rangle,
\end{equation}
where $\xi$ is the unit vector tangent to $\mathbb{R}$, and proved that its $(2,0)$-part (which for $n=2$ is just the
Abresch-Rosenberg differential) is holomorphic.

Using this quadratic form, we will prove the following.

\begin{theorem}\label{thm} Let $x:\Sigma^2\rightarrow M^n(c)\times\mathbb{R}$, $c\neq 0$, be an
isometrically immersed complete non-minimal pmc surface with
non-negative Gaussian curvature. Then one of the following holds:
\begin{enumerate}
\item the surface is flat;

\item $\Sigma^2$ is a minimal surface of a totally umbilical
hypersurface of $M^n(c)$;

\item $\Sigma^2$ is a cmc surface in a $3$-dimensional totally
umbilical submanifold of $M^n(c)$;

\item the surface lies in $M^4(c)\times\mathbb{R}\subset\mathbb{R}^6$
(endowed with the Lorentz metric), and there exists a plane $P$
such that the level lines of the height function
$p\rightarrow\langle x(p),\xi\rangle$ are curves lying in planes
parallel to $P$.
\end{enumerate}
\end{theorem}

\begin{remark}The same result was obtained by H. Alencar, M. do Carmo and R. Tribuzy in the case when $c<0$ (Theorem $3$ in \cite{AdCT}).
\end{remark}

In order to prove Theorem \ref{thm} we will need the following Simons type equation obtained by S.-Y. Cheng and S.-T. Yau (equation $2.8$ in \cite{CY}), which generalizes some previous results in \cite{NS,JS,S}. Let $N$ be an $n$-dimensional Riemannian manifold, and consider a symmetric operator $S$ on $N$, that satisfies the Codazzi equation $(\nabla_XS)Y=(\nabla_YS)X$, where $\nabla$ is the Levi-Civita connection on the manifold. Then, we have
\begin{equation}\label{delta}
\frac{1}{2}\Delta|S|^2=|\nabla S|^2+\sum_{i=1}^{n}\lambda_i(\trace S)_{ii}+\frac{1}{2}\sum_{i,j=1}^{n}R_{ijij}(\lambda_i-\lambda_j)^2,
\end{equation}
where $\lambda_i$, $1\leq i\leq n$, are the eigenvalues of $S$, and $R_{ijkl}$ are the components of the Riemannian curvature of $N$.

\section{The proof of Theorem \ref{thm}}

Let us consider an operator $S$, defined on the surface $\Sigma^2$ by
\begin{equation}\label{eq:S}
S=2A_H-c\langle T,\cdot\rangle T+\Big(\frac{c}{2}|T|^2-2|H|^2\Big)\id,
\end{equation}
where $T$ is the component of $\xi$ tangent to the surface. When the ambient space is $3$-dimensional this operator was introduced in \cite{B}. We shall prove that $|S|^2$ is a bounded subharmonic function on the surface.

First, it is easy to see that
\begin{equation}\label{eq:SQ}
\langle SX,Y\rangle=Q(X,Y)-\frac{\trace Q}{2}\langle X,Y\rangle,
\end{equation}
where $Q$ is the quadratic form given by \eqref{eq:Q}, which implies that $S$ is symmetric and traceless. Another direct consequence of \eqref{eq:SQ} is the following

\begin{lemma}\label{lemma:SQ0}  The $(2,0)$-part of $Q$ vanishes on $\Sigma^2$ if and
only if $S=0$ on the surface.
\end{lemma}

The following Lemma is proved in \cite{B}.

\begin{lemma}\label{lemma:Codazzi} The operator $S$ satisfies the Codazzi equation $(\nabla_XS)Y=(\nabla_YS)X$, where $\nabla$ is the Levi-Civita connection on the surface.
\end{lemma}

From Lemma \ref{lemma:Codazzi}, equation \eqref{delta} and the fact that $\trace S=0$, we easily get
\begin{equation}\label{eq:Simons}
\frac{1}{2}\Delta|S|^2=2K|S|^2+|\nabla S|^2,
\end{equation}
where $K$ is the Gaussian curvature of the surface.

Now, let us consider the local orthonormal frame field $\{E_3=\frac{H}{|H|},E_4,\ldots,E_{n+1}\}$ in the normal bundle, and denote $A_{\alpha}=A_{E_{\alpha}}$. It follows that $\trace A_3=2|H|$ and $\trace A_{\alpha}=0$, for all $\alpha>3$.

From the definition \eqref{eq:S} of $S$, we have, after a straightforward computation,
$$
\det A_3=\frac{1}{|H|^2}\det A_H=|H|^2-\frac{1}{8|H|^2}|S|^2-\frac{c^2}{16|H|^2}|T|^4-\frac{c}{4|H|^2}\langle ST,T\rangle,
$$
and then, by using the equation of Gauss of $\Sigma^2$ in $\bar M$,
$$
\begin{array}{ll}
R(X,Y)Z=&c\{\langle Y, Z\rangle X-\langle X, Z\rangle Y-\langle Y,T\rangle\langle Z,T\rangle X+\langle X,T\rangle\langle Z,T\rangle Y\\ \\&+\langle X,Z\rangle\langle Y,T\rangle T-\langle Y,Z\rangle\langle X,T\rangle T\}\\ \\&+\sum_{\alpha=3}^{n+1}\{\langle A_{\alpha}Y,Z\rangle A_{\alpha}X-\langle A_{\alpha}X,Z\rangle A_{\alpha}Y\}.
\end{array}
$$
The Gaussian curvature can be written as
\begin{equation}\label{K}
K=c(1-|T|^2)+|H|^2-\frac{1}{8|H|^2}|S|^2-\frac{c^2}{16|H|^2}|T|^4-\frac{c}{4|H|^2}\langle ST,T\rangle+\sum_{\alpha>3}\det A_{\alpha}.
\end{equation}

Since $\trace A_{\alpha}=0$, it follows that $\det A_{\alpha}\leq 0$, for all $\alpha>3$. Therefore, as $K\geq 0$, we get
$$
-\frac{1}{8|H|^2}|S|^2-\frac{c}{4|H|^2}\langle ST,T\rangle-\frac{c^2}{16|H|^2}|T|^4+c(1-|T|^2)+|H|^2\geq 0.
$$
From $|\langle ST,T\rangle|\leq\frac{1}{\sqrt{2}}|T||S|$ it results that $-\frac{c}{4|H|^2}\langle ST,T\rangle\leq\frac{|c|}{4\sqrt{2}|H|^2}|S|$, which implies
$$
-\frac{1}{8|H|^2}|S|^2+\frac{|c|}{4\sqrt{2}|H|^2}|S|+c(1-|T|^2)+|H|^2\geq 0.
$$

Next, we shall consider two cases as $c<0$ or $c>0$, and will prove that, in both situations, $|S|$ is bounded from above.

If $c<0$ we have
$$
-\frac{1}{8|H|^2}|S|^2-\frac{c}{4\sqrt{2}|H|^2}|S|+|H|^2\geq 0
$$
and then $|S|\leq \frac{\sqrt{c^2+|H|^2}-c}{\sqrt{2}}$.

When $c>0$ it follows that
$$
-\frac{1}{8|H|^2}|S|^2+\frac{c}{4\sqrt{2}|H|^2}|S|+c+|H|^2\geq 0,
$$
which is equivalent to $|S|\leq \frac{\sqrt{c^2+16c|H|^2+16|H|^2}+c}{\sqrt{2}}$.

As the surface is complete and has non-negative Gaussian curvature, it follows, from a result of A. Huber in \cite{H}, that $\Sigma^2$ is a parabolic space. From the above calculation and \eqref{eq:Simons}, we get that $|S|^2$ is a bounded subharmonic function and it follows that $|S|$ is a constant. Again using equation \ref{eq:Simons}, it concludes that $K=0$ or $S=0$. From Lemma \ref{lemma:SQ0}, we see that, when $\Sigma^2$ is not flat, the $(2,0)$-part of the quadratic form $Q$ vanishes on the surface, and then we obtain the last three items of our Theorem exactly as in the proofs of Theorem $2$ and Theorem $3$ in \cite{AdCT}.

\begin{remark} M. Batista characterized some cmc surfaces in
$M^2(c)\times\mathbb{R}$, under some assumptions on their mean
curvature and on $|S|$. Since these assumptions imply that
these surfaces have non-negative Gaussian curvature (this can be easily verified by using \eqref{K} and the fact that $|ST|^2=\frac{1}{2}|T|^2|S|^2$. We remark the
converse is not necessarily true), we can see that Theorem $3.1$ in \cite{ER}
generalizes his results (Theorem $1.2$ and Theorem $1.3$ in \cite{B}).
\end{remark}


\begin{thebibliography}{99}

\bibitem{AR} U. Abresch and H. Rosenberg, \textit{A Hopf differential for constant mean curvature
surfaces in $\mathbb{S}^2\times\mathbb{R}$ and
$\mathbb{H}^2\times\mathbb{R}$}, Acta Math. 193(2004), 141-–174.

\bibitem{AR2} U. Abresch and H. Rosenberg, \textit{Generalized Hopf
differentials}, Mat. Contemp. 28(2005), 1--28.

\bibitem {AdCT1} H. Alencar, M. do Carmo and R. Tribuzy, \textit{A
theorem of Hopf and the Cauchy-Riemann inequality}, Comm. Anal.
Geom. 15(2007), 283--298.

\bibitem{AdCT} H.~Alencar, M.~do Carmo and R.~Tribuzy, \textit{A
Hopf theorem for ambient spaces of dimensions higher than three},
J. Differential Geometry 84(2010), 1--17.

\bibitem{B} M.~Batista, \textit{Simons type equation in $\mathbb{S}^2\times\mathbb{R}$ and
$\mathbb{H}^2\times\mathbb{R}$ and applications}, preprint 2010.

\bibitem{CY} S.-Y. Cheng and S.-T. Yau, \textit{Hypersurfaces with constant scalar curvature}, Math. Ann. 225(1977), 195--204.

\bibitem{ER} J. M.~Espinar and H.~Rosenberg, \textit{Complete constant mean curvature
surfaces in homogeneous spaces}, Comment. Math. Helv., to appear.

\bibitem{H} A.~Huber, \textit{On subharmonic functions and differential geometry in the
large}, Comm. Math. Helv. 32(1957), 13--71.

\bibitem{NS} K.~Nomizu and B.~Smyth, \textit{A formula of Simons' type and hypersurfaces with constant mean
curvature}, J. Differential Geometry 3(1969), 367--377.

\bibitem{JS} J.~Simons, \textit{Minimal varieties in Riemannian
manifolds}, Ann. of Math. 88(1968), 62--105.

\bibitem{S} B.~Smyth, \textit{Submanifolds of constant mean curvature}, Math. Ann. 205(1973), 265--280.

\end{thebibliography}
\end{document}